\documentclass[a4paper,11pt]{article}

\usepackage[margin=3cm]{geometry}
\usepackage{fullpage}
\usepackage[english]{babel}
\usepackage[dvips]{epsfig}
\usepackage{amsfonts}
\usepackage{subfigure}
\usepackage{psfrag}
\usepackage{color}
\usepackage{soul}
\usepackage{theorem}

\usepackage[left]{lineno}

{\theoremstyle{plain}
\newtheorem{theorem}{Theorem}[section]
\newtheorem{proposition}[theorem]{Proposition}
\newtheorem{corollary}[theorem]{Corollary}
\newtheorem{lemma}[theorem]{Lemma}

{\theorembodyfont {\rmfamily}
\newtheorem{example}[theorem]{Example}
\newtheorem{nota}[theorem]{Notation}
\newtheorem{defn}[theorem]{Definition}
\newtheorem{remark}[theorem]{Remark}
\newtheorem{ex}[theorem]{Example}
}}

\begin {document}


\title{On the superimposition of Christoffel words}
\author{Genevi\`eve Paquin\thanks{with the support of FQRNT (Qu\'ebec)} \thanks{Laboratoire de math\'ematiques, CNRS UMR 5127, Universit\'e de Savoie, 73376 Le Bourget-du-lac Cedex, FRANCE, Genevieve.Paquin@univ-savoie.fr}
\and Christophe Reutenauer \thanks{Laboratoire de Combinatoire et d'Informatique Math\'ematique, Universit\'e du Qu\'ebec \`a Montr\'eal, CP. 8888 Succ. Centre-Ville, Montr\'eal, (QC) CANADA, H3C 3P8, Reutenauer.Christophe@uqam.ca } }

\maketitle

\newcommand{\QED}{\rule{1ex}{1ex} \par\medskip}

\newcommand{\nl}{\par\medskip\noindent}

\newcommand{\Proof }{{\nl\it Proof.\ }}

\newcommand{\fr}{\textnormal{Fr}}

\newcommand{\Card}{\hbox{\rm Card }}

\newcommand{\N}{\mathbb N}
\newcommand{\Q}{\mathbb Q}
\newcommand{\R}{\mathbb R}
\newcommand{\Z}{\mathbb Z}
\newcommand{\A}{\mathcal{A}}

\newcommand{\mod}{\textnormal{mod}}

\begin{abstract}  Initially {stated} in terms of Beatty sequences, the Fraenkel conjecture can be reformulated as follows: for a $k$-letter alphabet $\A$, with a fixed $k \geq 3$, there exists a unique balanced infinite word, up to letter permutations and shifts, that has mutually distinct letter frequencies. Motivated by the Fraenkel conjecture, we study in this paper whether two Christoffel words can be superimposed. Following from previous works on this conjecture using Beatty sequences, we give a necessary and sufficient condition for the superimposition of two Christoffel words having same length, and more generally, of two arbitrary Christoffel words. Moreover, for any two superimposable Christoffel words, we give the number of different possible superimpositions and we prove that there exists a superimposition that works for any two superimposable Christoffel words. Finally, some new properties of Christoffel words are obtained as well as a geometric proof of a classic result concerning the money problem, using Christoffel words. 
\end{abstract}


\section{Introduction}

{\it Beatty sequences} {and} {\it Sturmian words} {are equivalent objects. The first ones are studied in number theory. The second ones, known since the work of Bernoulli} \cite{jb1772}, {are studied in combinatorics on words and related domains.}  A {\it Beatty sequence} is a sequence of the form $S(\alpha, \beta)=\{\lfloor \alpha n + \beta \rfloor: n \in \Z\}$, with $\alpha, \beta \in \R$. It appeared in the literature for the first time in \cite{sb1926}, and the name came only more than 30 years later in \cite{igc1959, igc1960}. A finite set of Beatty sequences is called an {\it (eventual) exact cover} if every (sufficiently large) positive integer occurs in exactly one Beatty sequence. It is thus natural to wonder which sets of Beatty sequences are an (eventual) exact cover of the integers. Some particular cases have been studied for instance in \cite{sb1926, jvu1927, tb1957, ts1957,  rlg1963,in1963, af1969, fls1972}. Later, in \cite{eg1980},  appears the Fraenkel conjecture in terms of Beatty sequences which states that if a finite set of {\it rational Beatty sequences}, that is Beatty sequences with $\alpha \in \Q$,  is an eventual cover of the integers, then the $\alpha$'s satisfy a particular form (see \cite{eg1980} for more details). 

In combinatorics on words, the conjecture can be restated as: for a finite $k$-letter alphabet, with a fixed $k\geq 3$, there exists a unique balanced infinite word, up to letter permutations and shifts, that has mutually distinct letter frequencies. This supposedly unique infinite word is called a {\it Fraenkel word} and is given by the periodic word $\fr_k^\omega$, where $\fr_k$ is defined recursively by $\fr_k= \fr_{k-1}k\fr_{k-1}$ for $k\geq 2$, and $\fr_1=1$. 

Particular cases of the Fraenkel conjecture have been well studied  for instance by Morikawa, who published a series of papers on the topic  (see \cite{rt2001} for a good survey). More precisely, in \cite{rm1985}, the author proves the following theorem, which is a necessary and sufficient condition for the disjointness of two Beatty sequences, and that was later reformulated by Simpson \cite{rjs2004} as:

{\bf Theorem.} [\cite{rm1985}, {see also} \cite{rjs2004}] {\it Let $p_1,p_2,q_1,q_2$ be integers and let $p=\gcd(p_1,p_2)$, $q=\gcd(q_1,q_2)$, $u_1=q_1/q$ and $u_2=q_2/q$. There exist $\beta_1$ and $\beta_2$ such that the Beatty sequences $S_1=\{\lfloor p_1n/q_1+\beta_1\rfloor:n\in Z\}$ and $S_2=\{\lfloor p_2n/q_2+\beta_2\rfloor:n\in Z\}$ are disjoint if and only if there exist positive integers $x$ and $y$ such that
$$xu_1+yu_2=p-2u_1u_2(q-1).$$
}
This result is a step towards Fraenkel conjecture. In \cite{rjs2004}, {Simpson works out the proof of Morikawa, gives a new proof} and proves some new intermediate results. While translating Simpson results in terms of Christoffel words, some nice properties of these words naturally appear. In our paper, we first introduce some basic definitions and notation, and we show how the Fraenkel conjecture and the superimposition of  Christoffel words are related. Then after having formulated and proved the main results of Simpson in terms of Christoffel words, we go farther and give the number of superimpositions of two Christoffel words and one possible shift needed to superimpose them. We end this paper by showing how the geometric representation of Christoffel words can be used to prove a problem related to the {classical} money problem.

{The authors are grateful to R. J. Simpson for giving them useful details about his proofs. }

\section{Preliminaries}

We first recall notions on words (for more details, see for instance \cite{ml2002}).

An {\it alphabet} $\A$ is a finite set of symbols called {\it letters}. A {\it word over $\A$} is a sequence of letters from $\A$. The {\it empty word}  $\varepsilon$ is the empty sequence. Equipped with the concatenation operation, the set $\A^*$ of {\it finite words} over $\A$ is a free monoid with neutral element $\varepsilon$ and set of generators $\A$, and $\A^+=\A^* \setminus \varepsilon$. Given a nonempty finite word $u=u[0]u[1]\cdots u[n-1]$, with $u[i] \in \A$, the {\it length} $|u|$ of $u$ is the integer $n$. One has $|\varepsilon|=0$. We denote by $\A^n$ the set of finite words of length $n$ over $\A$ and by $\A^\omega$ the set of {\it (right-) infinite words} over $\A$.  The set $\A^\infty$ is defined as the set of finite and infinite words: $\A^\infty= \A^* \cup \A^\omega$. 

As usual, for a finite word $u$ and a positive integer $n$, the {\it $n$th power of $u$}, denoted $u^n$, is the word $\varepsilon$ if $n=0$; otherwise $u^n=u^{n-1}u$. The finite word $w$ is {\it primitive} if it is not the power of a shorter word. If $u\neq \varepsilon$, $u^\omega$ (resp. $\,^\omega u$) denotes the right-infinite (resp. left-infinite) word obtained by infinitely repeating $u$ to the right (resp. to the left). A right-infinite  word $\bf u$ is {\it periodic} (resp. {\it ultimately periodic}) if it can be written as ${\bf u}=w^\omega$ (resp. ${\bf u}=vw^\omega$), with $v\in \A^*$ and $w\in \A^+$. The set of {\it bi-infinite words}, denoted by $\A^\Z$, is defined as the set of functions $\Z \rightarrow \A$. For the sake of clarity, we denote in bold character a letter denoting an infinite or bi-infinite word, in opposition to a finite word. If $u \in \A^*$, then $\, ^\omega u ^\omega$ is the bi-infinite word ${\bf s}=\cdots u \bullet uu\cdots$. The point $\bullet$  is located between ${\bf s}[-1]$ and ${\bf s}[0]$ and represents the {\it origin} of the word $\bf s$. 

The number of occurrences of the letter $a$ in the word $u$ is denoted by $|u|_a$. Over infinite {words}, the {\it shift operator} $\sigma$ is defined by $\sigma: \A^\N \rightarrow \A^\N$ such that $\sigma({\bf s}[n])={\bf s}[n+1]$. It is also naturally defined over the set of bi-infinite words by $\sigma: \A^\Z \rightarrow \A^\Z$, with $\sigma({\bf s}[n])={\bf s}[n+1]$. A shift $\sigma ^k$, with $k\geq 0$, over a bi-infinite word is equivalent to move the origin $k$ times to the right. 

If, for some words $u, s \in \A^\infty$, $v, p \in \A^*$, $u=pvs$, then $v$ is a {\it factor} of $u$, $p$ is a {\it prefix} of $u$ and $s$ is a {\it suffix} of $u$. If $v\neq u$ (resp. $p\neq u$ and $s \neq u$), $v$ is called a {\it proper factor} (resp. {\it proper prefix} and {\it proper suffix}). The {\it set of factors} of the word $u$ is denoted $F(u)$.  

The {\it reversal} of the finite word $u=u[0]u[1]\cdots u[n-1]$, also called the {\it mirror image}, is $\tilde u= u[n-1]u[n-2]\cdots u[0]$ and if $u=\tilde u$, then $u$ is called a {\it palindrome}. Let $u \in \A^n$. Then $u[i]u[i+1]\cdots u[n-1]u[0]\cdots u[i-1]$ is a {\it conjugate} of $u$, for all $0\leq i \leq n-1$. 

In what follows, for $p,q \in \N$, we write $p \perp q$ if $\gcd(p,q)=1$. Otherwise, we write $p \not \perp q$. 

\subsection{Christoffel words}

In combinatorics on words, instead of using Beatty sequences, we use an equivalent combinatoric object: the Sturmian words. There exists a wide literature about Sturmian words in which we can find several characterizations depending on the context of the study {(see for instance} \cite{ml2002}). In particular, the Sturmian words are known as the balanced non-periodic infinite words over a $2$-letter alphabet. Recall that a {finite or infinite} word $w$ is {\it balanced} if for all finite factors $u,v \in F(w)$ having same length and for  all letters $a \in \A$, $\big | |u|_a - |v|_a\big | \leq 1$.  

A finite version of the Sturmian words is the family of Christoffel words. It has been studied for instance in \cite{ebc1875,br2006,bdlr2007, kr2007}. From the definition of Christoffel words given in \cite{mh1940}, in terms of symbolic dynamics, one can easily deduce the following:

\begin{defn} \label{defChristo} Let $\A=\{a<x\}$, $\alpha, \beta \in \N$ such that $\alpha \perp \beta$ and let  $n=\alpha + \beta$. The {\it  Christoffel word} $u\in \A^*$ with $\alpha$ occurrences of $a$'s and $\beta$ occurrences of $x$'s is defined by  $u=u[0]u[1]\cdots u[n-1]$, where 
$${\displaystyle u[i]= \left \{ \begin{array}{
ll} a & \textnormal{if }  (i+1)\beta \, \,  \mod \, \, n > i\beta \, \, \mod \, \,  n\\
x & \textnormal{if } (i+1)\beta \, \, \mod \, \, n \leq i\beta \, \,  \mod\, \,  n
\end{array} 
      \right . }$$
for $0 \leq i < n$, where $i\beta \, \, \mod \, \,  n$ denotes the {remainder} of the Euclidean division of $i\beta$ by $n$. We say that $u$ has {\it slope} $\beta/\alpha$. 
\end{defn}

To any Christoffel word, we can associate a Christoffel path, defined as follows.  

\begin{defn} Suppose $a,b \in \N$ and $a \perp b$. The {\it Christoffel path} of {\it slope} $b/a$ is the path from $(0,0)$ to $(a,b)$ in the integer lattice $\Z \times \Z$ that satisfies the following two conditions. 
\begin{itemize} 
\item [i)] The path lies below the line segment that begins at the origin and ends at $(a,b)$.
\item [ii)] The region in the plane enclosed by the path and the line segment contains no other points of $\Z \times \Z$ besides those of the path.
\end{itemize}
\end{defn}

The next figure shows the Christoffel path of slope $3/5$. 
\begin{center}
\includegraphics[width=5cm]{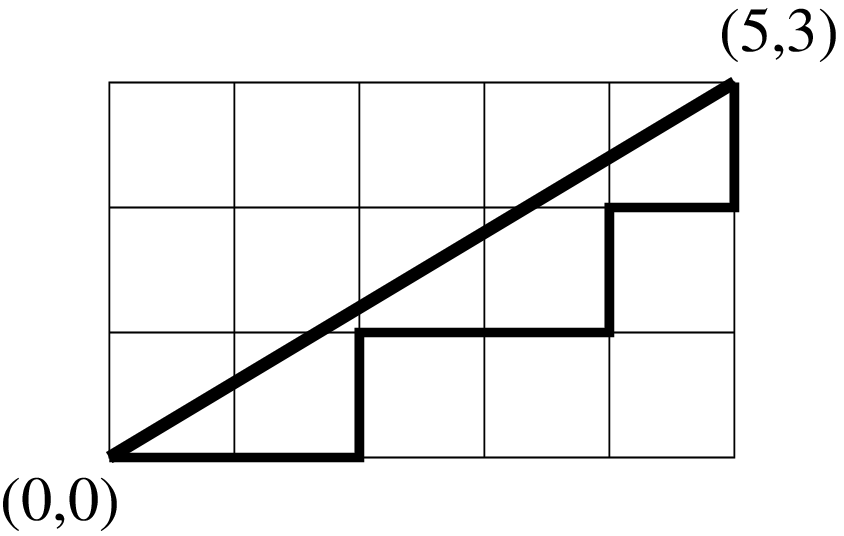}
\end{center}

Notice that Definition \ref{defChristo} can be generalized to powers of Christoffel words by removing the condition $\alpha \perp \beta$. We then have:

\begin{defn}  Let $C(n,\alpha)$ be a word of length $n$ over $\{a<x\}^*$ having $\alpha$ occurrences of $a$'s, {with $\alpha \leq n$,} and let $r=\gcd(n,\alpha)$.
\begin{itemize}
\item [i)] If $r=1$, then $C(n,\alpha)$ denotes the Christoffel word of slope $\frac{n-\alpha}{\alpha}$.
\item [ii)] If $r >1$, then $C(n,\alpha)=C(r\frac{n}{r},r\frac{\alpha}{r})=\left(C\left (\frac{n}{r},\frac{\alpha}{r}\right )\right)^r$ denotes the $r$th power of the Christoffel word  $C\left (\frac{n}{r},\frac{\alpha}{r}\right )$.
\end{itemize}
\end{defn}

{Let us recall the following classical lemma.}
\begin{lemma}  \label{cayleyMir} The reversal of a Christoffel word (resp. power of a Christoffel word) $C(n,\alpha) \in \{a< x\}^*$, denoted by $\widetilde C(n,\alpha)$, is also a Christoffel word (resp.  power of a Christoffel word) over the same alphabet, but for which the order of the letters is reversed. More precisely,  $\widetilde C(n,\alpha)=C(n,n-\alpha) \in \{x< a\}^*$.
\end{lemma}

Let us now consider the directed graph with the set of vertices $\{0,1,2, \ldots, \alpha +\beta-1\}$ that has an arrow from the vertex $i$ to the vertex $j$ if $i+\beta \equiv j \, \mod \, n$ and labeled by $a$ if $i < j$, and by  $x$ if $j< i$.\begin{itemize}
\item [\rm i)] If $\alpha \perp \beta$, this graph is called the {\it  Cayley graph} of the Christoffel word $u$  over the alphabet $\{a<x\}$ and having slope $\displaystyle \frac{\beta}{\alpha}$. 
\item [\rm ii)] If $\gcd(\alpha, \beta)=r>1$, then the graph obtained is isomorphic to the Cayley graph $\displaystyle C\left (\frac{\alpha+\beta}{r},\frac{\alpha}{r}\right )$. This graph read $r$ times is the Cayley graph of $C(\alpha+\beta, \alpha)$.
\end{itemize}

\begin{example} The Cayley graph associated to the Christoffel word over $\{a< x\}$ having slope $3/5$ is
\begin{center}
\includegraphics[width=5cm]{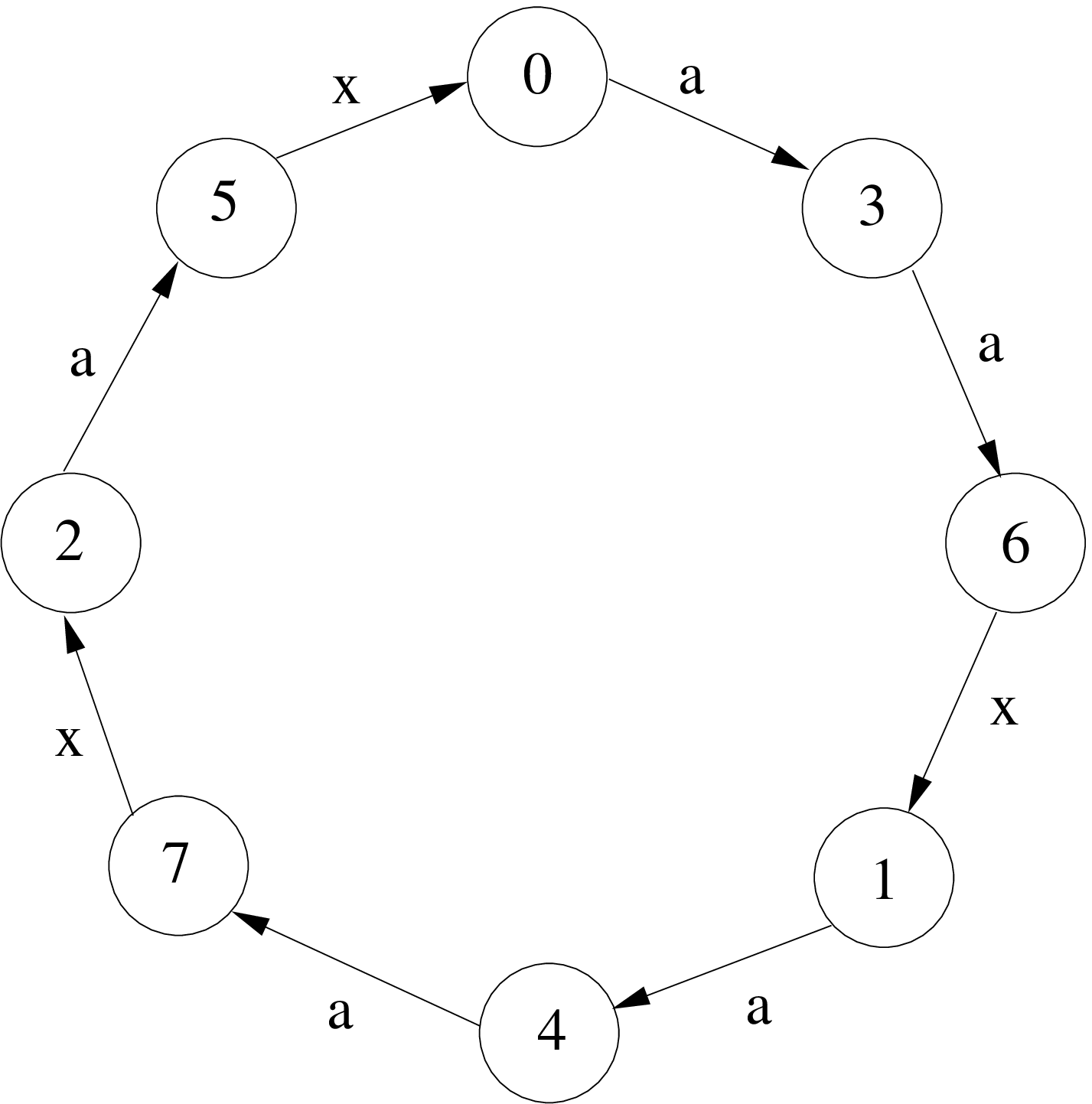}
\end{center}
and  $u=aaxaaxax$, $|u|_a=5, |u|_x=3$.
\end{example}

\subsection{Link with the Fraenkel conjecture}

Before showing how the problem of superimposition of Christoffel words is related to the Fraenkel conjecture, some definitions are required. 

First, let us recall that a word  $w \in \A^*$ is said to be {\it circularly balanced} if  $w^2=ww$ is balanced.  

{\ex The word $u=112121$ is balanced, but is not circularly balanced. Indeed, $111, 212 \in F(uu)$, but $\big | |111|_1 - |212|_1\big | >1$.
}

{\ex One can easily verify that the word $112112$ is balanced and consequently, that $v=112$ is circularly balanced. 
}

Let $w \in \A^*$, with $\Card \A \geq 3$. The  {\it projection}  $\Pi_a(w)$ of the word $w \in \A^*$ {on} the alphabet $\{a,x\}$, with $a \in \A$ and $x \notin \A$, is defined by
$${\displaystyle \Pi_a (w)[i]= \left \{ \begin{array}{
ll} a & \textnormal{if }  w[i]=a\\
x & \textnormal{otherwise.}
\end{array} 
      \right . }$$

{\ex Let $w=1232343112$. Then $\Pi_1(w)=1xxxxxx11x$, $\Pi_2(w)=x2x2xxxxx2$, $\Pi_3(w)=xx3x3x3xxx$ and $\Pi_4(w)=xxxxx4xxxx$. 
}

{The next result is given without proof, since it is trivial.}
\begin{lemma} \label{proj} If $w \in \A^*$ is circularly balanced, then for all $a \in \A$, the projection $\Pi_a(w)$ is so. 
\end{lemma}

\begin{defn} {\bf (superimposition of bi-infinite words)} Let ${\bf u} \in \,^\omega\{a,x\}^\omega$ and ${\bf v} \in \,^\omega\{b,x\}^\omega$ be two bi-infinite words. Let $A$ be the set of positions of the $a$'s in ${\bf u}$ and $B$ be the set of positions of the $b$'s in ${\bf v}$. We say that ${\bf u}$ and ${\bf v}$ are {\it superimposable} if $A \cap B =\emptyset$.
\end{defn}

{\ex Let ${\bf u}= \hspace{-0.1cm}\,^\omega(aaxaxx)^\omega$ and ${\bf v}=\hspace{-0.1cm}\,^\omega(xxbxxx)^\omega$. Then $A=\{0,1,3\}+6\Z$ and $B=\{2\} +6\Z$. Hence ${\bf u}$ and ${\bf v}$ are superimposable. 
}

\begin{defn} \label{defSuperp} {\bf (superimposition of finite words)} Let $u \in \{a,x\}^n$ and $v \in \{b,x\}^m$ be finite words. Let $A$ be the set of positions of the  $a$'s in  $u$ and let $B$ be the set of  positions of the $b$'s in $v$.  Then $u$ and $v$ are {\it superimposable} if and only if there exists  $k \in \Z$ such that  $^\omega u ^\omega$ and $\sigma^k(^\omega v ^\omega)$ are  superimposable, that means if 
$$(A+n\Z) \cap ({B-k}+m\Z)=\emptyset.$$ 
If  $k=0$, $u$ and $v$ are said {\it perfectly superimposable}.
\end{defn}

\begin{remark} In Definition \ref{defSuperp}, the condition  $k \in \Z$ can be replaced by  $k \in [0, \min\{m,n\}-1]$. Indeed, one can easily verify that if there exists a shift $k$ outside this interval that allows the  superimposition, then there exists $k' \in [0,\min\{m,n\}-1]$ such that it is so. 
\end{remark}



For a letter $\alpha \in \A$ and a finite word $w\in \A^*$, the {\it conjugacy operator} $\gamma$ is defined by $\gamma(\alpha w)=w\alpha $. Then  $\sigma( w^\omega)=(\gamma w )^\omega$: $\gamma$ acts over finite words as the shift $\sigma$ acts over infinite words.

\begin{lemma} Let $u, v,  A$ and $B$ be such as in Definition \ref{defSuperp}. The words $u$ and $v$  are superimposable and are such that 
$$(A+n\Z)\cap (B-k+m\Z) =\emptyset$$
if and only if  $u$ and $\gamma^{k}(v)$ are perfectly superimposable.
\end{lemma}

\Proof By definition,  $u$ and $v$ are superimposable if and only if there exists $k \in \Z$ such that $(A +n\Z) \cap (B-k+m\Z)=\emptyset$. This condition is satisfied if and only if the set of positions of the $a$'s in $\,^\omega u^\omega $ and the set of positions of the $b$'s in $\sigma^k(\,^\omega v^\omega)$ are disjoint. It is then sufficient to show that the positions of the $b$'s in $\,^\omega( \gamma^{k}(v))^\omega$ are the same as in $\sigma^k(\,^\omega v ^\omega)$. Those two words have period $m$ and consequently, the positions of the $b$'s are the same if and only if $\, ^\omega(\gamma^k(v))^\omega=\sigma ^k(\,^\omega v^\omega) \iff \gamma^k(v)=\sigma^k(\, ^\omega v^\omega)[0, m-1]$. This last condition is satisfied by the definition of $\gamma$. \QED


\begin{corollary} \label{celui} A finite circularly balanced word $w \in \A^m$, with $\A=\{1,2,\ldots, k\}$,  having  pairwise distinct letter frequencies can be obtained by the superimposition of  $k$ circularly balanced words $w_1\in \{1,x\}^m, w_2 \in \{2,x\}^m, \ldots, w_k \in \{k,x\}^m$ such that  $|w_i|_i = |w|_i$ for all $1\leq i \leq k$.\end{corollary}

\Proof It is sufficient to apply the projection $\Pi_a(w)$ to all letters $a \in \A$ and to conclude using  Lemma \ref{proj}. \QED

\begin{lemma} {\rm \cite{agh1998}}\label{periodic} Any balanced infinite word over a $k$-letter alphabet, with $k\geq 3$, having pairwise distinct letter frequencies is periodic. 
\end{lemma}

Corollary \ref{celui} and Lemma \ref{periodic} give the main motivation of this paper. Lemma \ref{periodic} tells us that in order to prove the Fraenkel conjecture, it is sufficient to prove that for a $k$-letter alphabet, any circularly balanced finite word having pairwise distinct letter frequencies is a conjugate of $\fr_k$. Moreover, we deduce from Corollary \ref{celui}  that any finite word satisfying the conditions of the Fraenkel conjecture can be obtained  by the superimposition of circularly balanced words, or in other words, by the superimposition of conjugates of powers of Christoffel words, since Christoffel words are primitive balanced {words} that are minimal with respect to the lexicographic order in their conjugacy class.  


In this paper, we are naturally interested in the superimposition of two circularly balanced words. We  first  only consider two  finite primitive words and we give a necessary and sufficient condition for the superimposition of those two words. Moreover, if  $u$ is primitive and circularly balanced, then there exists a conjugate of  $u$ that is a Christoffel word. Thus, we consider the corresponding Christoffel words  $w_1$ and  $w_2$ and we give a criteria such that  $w_1$ and $w_2$ are superimposable. To do so, we use results from  \cite{rjs2004} which are an extension of the works of  \cite{rm1985}. We will see that considering a finite circularly balanced word as a conjugate of a power of a Christoffel word allows us to get some nice properties of the Christoffel words.

\section{Superimposition of Christoffel words having same length}

In this section, we first recall some properties of Christoffel words that will be {used} in the sequel. Then  we study the superimposition of Christoffel words having same length. Notice that most of the results of this section are already known from \cite{rm1985, rjs2004}, but we include some of their proofs since {it} will be  useful in order to prove the new results presented in the last subsection (number of superimpositions) and in the next sections. 

{In the sequel,  for a positive integer $\alpha$ and a fixed $n$, we will denote by $\overline \alpha$ the integer in $[0\ldots n-1]$ such that  $\alpha \overline \alpha \equiv -1 \, \mod \, n$.}

Lemma \ref{th3Simpson} is a translation of Theorem 3 in \cite{rjs2004} in terms of Christoffel words. This result also appears in  \cite{bdlr2007} in an equivalent form using the duality of Christoffel words.

\begin{lemma} {\rm \cite{rjs2004,bdlr2007} } \label{th3Simpson} Let $C(n,\alpha) \in \{ a < x\}^*$ be a Christoffel word. Then  the positions  of the  $a$'s  modulo $n$ in  $C(n,\alpha)$ are given by the set $\{0, \overline \alpha, 2\overline {\alpha}, \ldots, (\alpha-1)\overline \alpha \}$.
\end{lemma}


Lemma \ref{th3Simpson} can be easily generalized to a power of a Christoffel word as:

\begin{corollary} \label{corPos} Let $C(nq,\alpha q)=(C(n,\alpha))^q$ with $n \perp \alpha$. Then the positions of the $a$'s modulo $nq$ in  $C(nq,\alpha q)$ are given by 
$$\bigcup_{i=0}^{q-1}\{0,\overline \alpha, 2\overline \alpha, \ldots, (\alpha-1)\overline \alpha\} +in.$$
\end{corollary}

{The following theorem is deduced by Simpson from the Chinese remainder Theorem.}
 
\begin{theorem} {\rm \cite{rjs2004}} Let $C(n, 1)$ and $C(m, 1)$ be two Christoffel words. Then $C(n,1)$ and $C(m,1)$ are superimposable if and only if  $n\not \perp m$.
\end{theorem}


\begin{lemma} \label{lemIneq} Let  $C(n,\alpha)\in \{a<x\}^*$. For each position  $i$ of an $a$ in $\widetilde C(n,\alpha)$, there exists  $j \in \N$ such that 
$$i\alpha < jn \leq (i+1)\alpha.$$
\end{lemma}

\Proof Recall that by Lemma \ref{cayleyMir},  $\widetilde C(n,\alpha)=C(n,n-\alpha) \in \{x< a\}^*$. Using the generalization of Definition \ref{defChristo} to powers of Christoffel words by replacing respectively  $a, x,\alpha$ and $\beta$ by $x, a, 
n-\alpha$ and $\alpha$, we obtain that $\widetilde C(n,\alpha)[i]=a$ if and only if  $(i+1)\alpha \, \, \mod \, \,n \leq i\alpha \, \, \mod \, \,n$.  But $(i+1)\alpha \, \, \mod \, \, n \leq i\alpha \, \, \mod \, \,  n$ if and only if there exists a multiple of $n$ between $i\alpha$ and $(i+1)\alpha$ inclusively. This condition is satisfied if and only if there exists $j$ such that $i\alpha < jn \leq (i+1) \alpha$.  
{\flushright \QED}

\begin{lemma} \label{lempos}  Let $C(n, \alpha)\in \{a<x\}^n$  and $C(n,\beta)\in \{b<x\}^n$  be Christoffel words or power of Christoffel words. If $\alpha | \beta$, then the set of positions of the  $a$'s in $C(n,\alpha)$ is a subset of the set of positions of the $b$'s in  $C(n,\beta)$.
\end{lemma}

\Proof Let us prove this statement for the reversed words. Since $\alpha | \beta$, we write $\beta=q\alpha$, $q \in \N$. Let $i$ be the position of an $a$ in $\widetilde C(n,\alpha)$. Then by Lemma \ref{lemIneq}, there exists $j \in \N$ such that  $i\alpha <jn  \leq (i+1)\alpha$. Multiplying both sides of the inequation by $q$ yields 
$$i(q\alpha) < (jq)n\leq (i+1)(q\alpha) \iff i\beta< (jq)n \leq (i+1)\beta,$$
with $jq\in \N$. Hence $i$ is also a position of a $b$ in $\widetilde C(n,\beta)$. \QED

\begin{theorem}  [Th. 2 in \cite{rjs2004}] \label{th2Simpson} Let $C(n, \alpha)\in \{a<x\}^n$, $C(m, \beta)\in~\{b<~x\}^m$ be two Christoffel words and let  $p=\gcd(m,n)$. Then $C(n,\alpha)$ and $\gamma^{k}C(m,\beta)$ are perfectly superimposable if and only if  $C(p, \alpha)$ and $\gamma^{k}C(p, \beta)$ are so. 
\end{theorem}

The proof of Theorem \ref{th2Simpson} in terms of Christoffel words can be found in \cite{gp2008}.

\begin{corollary} \label{cor1} If the Christoffel words $C(n,\alpha)$ and $C(m,\beta)$ are superimposable, then $m\not \perp n$ and $\alpha + \beta  \leq p$, with $p=\gcd(m,n)$.
\end{corollary}

\Proof By Theorem \ref{th2Simpson},   $C(n,\alpha)$ and $C(m,\beta)$ are superimposable if and only if $C(p,\alpha)$ and $C(p,\beta)$ are so. This implies that if  $C(n,\alpha)$ and $C(m,\beta)$ are superimposable, then $\alpha +\beta \leq p$. Since $\alpha, \beta >0$, we have $1< \alpha +\beta \leq p$,  and consequently, $m\not \perp n$.  \QED

In what follows, we will first consider only Christoffel words having same length, since Theorem \ref{th2Simpson} will then allow us {in Section 4 to generalize our results to words of any length.} 

\subsection{Particular case: if $\alpha | \beta$}

In this subsection, we study the superimposition of the Christoffel words $C(n,\alpha)$ and $C(n,\beta)$, {\it having same length},  with $\alpha | \beta$. We give a criteria that the shift must satisfy in order to allow the superimposition (Lemma \ref{lemSup}), and then, we show a necessary and sufficient condition for the superimposition of those Christoffel words (Corollary \ref{cor5Simpson}). We also exhibit a shift that will always allow the perfect superimposition of two Christoffel words (Corollary \ref{uneSuperp}). We end the subsection by showing how a Christoffel word can be {viewed}  as the superimposition of some Christoffel words  (Theorem \ref{th6Simpson}).  

\begin{lemma} \label{lemMiroir} Let $C(n,\alpha)\in \{a<x\}^n$ be a Christoffel word. Then ${C}(n,\alpha)=\gamma^ {\overline{\alpha}}\widetilde C(n,\alpha)$.
\end{lemma}

\Proof By Lemma \ref{th3Simpson}, the positions of the $a$'s (modulo $n$) in $C(n,\alpha)$ are given by the set $$A=\{0, \overline \alpha, 2\overline \alpha, \ldots, (\alpha-1)\overline \alpha\}.$$ In the other hand, the positions of the $a$'s (modulo $n$) in the reverse word $\widetilde C(n,\alpha)$ are given by 
\begin{eqnarray*}
\widetilde A&=&\{n-1, n-1-\overline \alpha, n-1-2\overline \alpha, \ldots, n-1-(\alpha-1)\overline \alpha\}\\
&=&\{-1, -1-\overline \alpha, -1-2\overline \alpha, \ldots, -1-(\alpha-1)\overline \alpha\}.
\end{eqnarray*}
We then obtain that the positions of the $a$'s in the conjugate $\gamma^{\overline \alpha}\widetilde C(n,\alpha)$ are given by
\begin{eqnarray*}
\gamma^{\overline \alpha}\widetilde A &=&\{-1-\overline \alpha, -1-2\overline \alpha, -1-3\overline \alpha , \ldots, -1-\alpha\overline \alpha \}\\
&=&\{(\alpha-1)\overline \alpha, (\alpha-2)\overline \alpha, \ldots,,\overline\alpha, 0\}=A.
\end{eqnarray*}
\vspace{-1.7cm}
{\flushright \QED}

\begin{defn} Let $I=[a,b]$ and $I'=[c,d]$ be two intervals of integers. We say that $I$ is located {\it at the left} of  $I'$ if $a<c$. 
\end{defn}

\begin{lemma} \label{inter} The set of differences between the positions of the $a$'s in $C(n,\alpha) \in \{a<x\}^n$ and the positions of the $b$'s in $C(n,\beta)\in \{b<x\}^n$, with $\beta=q\alpha$ and $q \in \N$, forms an interval of integers having cardinality $(2\alpha -1)q$. 
\end{lemma}

\Proof Recall from Lemma \ref{th3Simpson} that the  positions of the $a$'s in $C(n,\alpha)$ are $\{0,\overline \alpha,\ldots,(\alpha-1)\overline \alpha\}$ and those of the $b$'s in $C(n,\beta)$ are $\{0,\overline \beta,\ldots,(\beta-1)\overline \beta\}$. Since $\beta=q\alpha$, multiplying both sides by  $\overline \alpha \overline \beta$ yields  $\overline \alpha\equiv q\overline \beta \, \mod \, n $ and hence $i\overline \alpha\equiv iq\overline \beta \, \mod \, n$.  Consequently, the differences between the positions of the letters form the set
\begin{equation} \label{eq74}
 E=\{j\overline \beta-i\overline \alpha\}_
{0\leq j < \beta \atop 0 \leq i < \alpha}
=
\{(j-iq)\overline \beta\}_
{0\leq j < \beta \atop
0 \leq i < \alpha
}.
\end{equation}
For a fixed $i$, the possible values of $j-iq$ form the interval $[-iq, \beta -iq[$. Since $q>0$, we deduce that for any $i$, the interval $[-iq, \beta-iq[$ is at the left of the interval $[-(i-1)q, \beta-(i-1)q[$. 

We have $\beta=q\alpha \Rightarrow \beta \geq q \Rightarrow \beta -iq \geq  q-iq=-(i-1)q$. Thus, the union of two consecutive intervals is also an interval, and consequently, the union of these $\alpha$ intervals forms the interval $[-(\alpha-1)q, \beta[$, which has cardinality
\begin{equation}\label{eq741}
\beta-\left( -(\alpha-1)q\right )=\beta+\alpha q -q=\alpha q+\alpha q -q= (2\alpha -1)q.
\end{equation}
\vspace{-1cm}{\flushright \QED}

\begin{lemma} [Th. 4 in \cite{rjs2004}] \label{lemSup} Let $C(n,\alpha)\in \{a<x\}^n$ and $C(n,\beta)\in \{b<x\}^n$ be two  Christoffel words, with $\beta=q\alpha$ and $q \in \N$ and let $\ell \in [0,n-1]$. The following conditions are equivalent:
\begin{itemize}
\item [i)] $C(n,\alpha)$ and $\gamma^{\ell\overline \beta}C(n,\beta)$ are perfectly superimposable;
\item [ii)]  $\ell+n\N \cap [-(\alpha-1)q,\beta[  =\emptyset$; 
\item [iii)] $C(n,\alpha)$ and $\gamma ^{\overline \beta (1+\ell)}\widetilde C(n,\beta)$ are perfectly superimposable.
\end{itemize}
\end{lemma}

\Proof $C(n,\alpha)$ and $\gamma^{k}C(n,\beta)$ are perfectly superimposable if and only if the shift  $k$ is not contained in the set $E$ (see (\ref{eq74})). Otherwise, there is a $a$ in $C(n,\alpha)$ at the same position as a $b$ in $\gamma^k C(n,\beta)$. This last condition is satisfied  if and only if there exits $\ell \notin [-(\alpha-1)q,\beta[ \, \mod \, n$ such that $k=\ell \overline \beta$. Thus there exists  $\ell\notin [-(\alpha-1)q, \beta[ \, \mod \, n$ if and only if $C(n,\alpha)$ and $\gamma^{\ell \overline{\beta}}C(n,\beta)$ are perfectly superimposable. Hence i) $\iff$ ii). Moreover, Lemma \ref{lemMiroir} gives that $C(n,\beta)=\gamma^{\overline \beta}\widetilde{C}(n,\beta)$. Replacing  $C(n,\beta)$ by this value in  $\gamma ^{\ell \overline \beta}C(n,\beta)$ yields
$$ \gamma^{\ell \overline \beta}C(n,\beta)= \gamma^{\ell \overline \beta}\gamma^{\overline \beta}\widetilde C(n,\beta) = \gamma ^{\overline \beta (\ell +1)}\widetilde C(n,\beta).$$
Hence  i) $\iff$ iii). \QED

\begin{corollary} [Cor. 5 in \cite{rjs2004}] \label{cor5Simpson} Let $C(n, \alpha)\in \{a<x\}^n$ and $C(n, \beta) \in \{b<x\}^n$ be Christoffel words such that  $\beta=q\alpha$, $q \in \N$. Then $C(n,\alpha)$ and $C(n,\beta)$ are superimposable if and only if  $(2\alpha -1)q < n$.
\end{corollary}

\Proof  The words $C(n,\alpha)$ and $C(n,\beta)$ are superimposable if and only if there exists a shift  $0\leq k < n$ such that the positions of the $\alpha$ occurrences of  $a$'s in $C(n,\alpha)$ form a disjoint set from the set of positions of the $\beta$ occurrences of $b$'s in $C(n,\beta)$. Such a shift $k$ exists if and only if the set $E$ (from Equation (\ref{eq74})) has cardinality at most $n-1$. We conclude using the fact that by Equation (\ref{eq741}),  $\Card(E)=(2\alpha -1)q$. \QED

From Lemma \ref{lemSup}, it is also possible to deduce a shift that always allows the perfect superimposition of two superimposable Christoffel words $C(n,\alpha)$ and $C(n,\beta)$ having same length, with $\alpha | \beta$:

\begin{corollary} \label{uneSuperp} Let $C(n,\alpha)$ and $C(n,\beta)$ be two superimposable Christoffel words such that  $\beta=q\alpha $, $q \in \N$. Then $C(n,\alpha)$ and $\gamma^{(1-r)}\widetilde C(n,\beta)$ are perfectly superimposable, with $\alpha r \equiv 1 \, \mod \, n$.
\end{corollary}

\Proof By Lemma \ref{lemSup}, it is sufficient to prove that there exists $\ell \notin [- (\alpha-1)q, \beta[ \, \, \mod \,  n$ such that 
$$\overline \beta(1+\ell) \equiv 1- r\, \mod \, n.$$
Multiplying both sides by  $\beta$ yields
$$-1-\ell \equiv \beta -\beta r \, \mod \,  n \iff  \ell \equiv \beta r -\beta -1 \, \mod \,  n.$$
Modulo $n$, we have
$$\beta r -\beta -1 \notin [-(\alpha-1)q,\beta[ \,\,= [-\beta+q,\beta[ \iff  \beta r\notin [q+1,2\beta +1[.$$
Since $\beta=q\alpha$ and $\alpha r \equiv 1 \, \mod \,  n$, we have  $\beta r \equiv q \, \mod \,  n$.
It is then sufficient to show that 
$$q +n \geq 2\beta +1 \iff n \geq 2\alpha q+1-q=(2\alpha -1)q+1.$$
We conclude using Corollary \ref{cor5Simpson}. \QED

\begin{theorem}[Th. 6 in \cite{rjs2004}] \label{th6Simpson} Let $C(n,q\alpha)\in\{a<x\}^n$ be a Christoffel word. Then the set of positions of the $a$'s in $C(n,q\alpha)$ is the union of the sets  $\{0,\overline \alpha,...,(\alpha-1)\overline \alpha\}+k\overline{q\alpha},$
for $0 \leq k < q$. Moreover, the Christoffel word $C(n,q\alpha)$ is the result of the exact superimposition of the  following $q$ conjugates of $C(n,\alpha)$: $C(n,\alpha), \gamma^{-\overline {q\alpha}}C(n,\alpha),\ldots,\gamma^{-(q-1)\overline{q\alpha}}C(n,\alpha).$
\end{theorem}

\Proof By Lemma \ref{th3Simpson}, the set of positions of the $a$'s in $C(n,q\alpha)$ is, modulo $n$,  
\begin{equation}\{0,\overline{q\alpha}, 2\overline{q\alpha},\ldots, (q\alpha-1)\overline{q\alpha}\}= \bigcup_{j=0}^{q\alpha -1}j\overline{q\alpha}=\bigcup_{k=0}^{q-1}\bigcup_{i=0}^{\alpha-1} (iq+k)\overline{q \alpha}. \label{eqnunion}
\end{equation}
The last equality is obtained by separating the positions with respect to their remainder modulo $q$. Since $q \alpha \overline{q \alpha}\equiv -1 \, \mod \, n$, we have $q\overline {q \alpha}\equiv \overline \alpha \, \mod \, n$. Thus replacing $q\overline{q\alpha}$ by $\overline \alpha$ in Equation (\ref{eqnunion}) yields
\begin{eqnarray}
\bigcup_{j=0}^{q\alpha -1}j\overline{q\alpha}&=&\bigcup_{k=0}^{q-1}\bigcup_{i=0}^{\alpha-1} i\overline \alpha +k\overline{q\alpha}= \bigcup_{k=0}^{q-1}  \{0, \overline \alpha, 2\overline \alpha, \ldots, (\alpha-1)\overline \alpha \}+k\overline{q\alpha}.\label{eqnunion2}
\end{eqnarray}

We conclude this proof by observing that the $q$ sets $\{0, \overline \alpha, \ldots, (\alpha-1)\overline \alpha\} + k\overline{q\alpha}$, for $0\leq k \leq q-1$, 
 correspond respectively to the positions of the $a$'s in the conjugates of Christoffel words $C(n,\alpha),$ $\gamma^{-\overline{q\alpha}}C(n,\alpha),$ $\ldots, \gamma^{-(q-1)\overline{q\alpha}}C(n,\alpha)$.  \QED

\subsection{General case} \label{secmemelng}

In this section, we study the general case of the superimposition of two Christoffel words having same length. In order to do so,  we consider the Christoffel words $C(n,q\alpha)\in \{a< x\}^*$ and $C(n,q\beta)\in \{b< x\}^*$, with $\alpha \perp \beta$ and $q \in \N$. 

\begin{nota} For $0 \leq i < \alpha $, we denote by $V_i$ the interval of integers
$$V_i=[(-q+1)\beta, q\beta -1] +i\overline \alpha \beta. $$
\end{nota}

\begin{proposition} \label{unionInter} The Christoffel words $C(n,q\alpha)$ and $C(n,q\beta)$ are superimposable if and only if the union 
\begin{equation} \label{unionVi} \displaystyle \bigcup _{i=0}^{\alpha-1} V_i \end{equation}
is not a complete set of residues modulo $n$. 
\end{proposition}

\Proof By inverting $q$ and $\alpha$ in Theorem \ref{th6Simpson}, we find that $C(n,q\alpha)$ is the perfect superimposition of the $\alpha$ conjugates $C(n,q), \gamma^{-\overline{q\alpha}}C(n,q), \ldots, \gamma^{-(\alpha-1)\overline{q\alpha}}C(n,q)$. The set of  positions of the $a$'s in  $\displaystyle C(n,q\alpha)$ is $\bigcup_{i=0}^{\alpha-1}\textnormal{pos}_a(\gamma^{-i\overline{q\alpha}}C(n,q))$, where pos$_a(w)$ denotes the positions of the $a$'s in $w$. Moreover, replacing $\alpha$, $q$ and $\beta$ by respectively $q$, $\beta$ and $q\beta$ in Lemma \ref{lemSup} yields that   $C(n,q)$ and $\gamma^{\ell \overline{q\beta}}C(n,q\beta)$ are perfectly superimposable if and only if there exists  $\ell \notin [-(q-1)\beta, q\beta[ \,  \mod \,  n$. More generally, $\gamma^{-i{\overline{q\alpha}}}C(n,q)$ and $\gamma^{\ell \overline{q\beta}}C(n,q\beta)$ are perfectly superimposable if and only if $C(n,q)$ and $\gamma^{\ell \overline{q\beta}+i\overline{q\alpha}}C(n,q\beta)$ are perfectly superimposable. In order to get the form of Lemma \ref{lemSup} iii), we rewrite $\ell \overline{q\beta}+i\overline{q\alpha}$ as 
\begin{eqnarray*}
\ell \overline{q\beta}+i\overline{q\alpha}&=&\ell \overline{q\beta}-\overline{q\beta}q\beta i \overline{q\alpha}\\
&=& \overline{q\beta}(\ell +q\beta i \overline q \,  \overline \alpha)\\
&=& \overline {q\beta}(\ell -i\overline \alpha \beta).
\end{eqnarray*}
We now have the required form of Lemma \ref{lemSup} iii). Then    $\gamma^{\ell \overline{q\beta}+i\overline{q\alpha}}C(n,q\beta)$ and $C(n,q)$ are perfectly superimposable if and only if there exists $\ell -i\overline \alpha \beta \notin [-(q-1)\beta, q\beta[ \, \mod \,  n$. This last condition is equivalent to the existence of a $\ell \notin  [-(q-1)\beta, q\beta[ \, +i\overline \alpha \beta=V_i$, but we need that  $\ell \notin V_i$ for all $0 \leq i < \alpha$. Thus the words $C(n,q\alpha)$ and $C(n,q\beta)$ are superimposable if and only if $\bigcup_{i=0}^{\alpha-1}V_i$ is not a complete set of residues modulo $n$. \QED

\begin{corollary} \label{corTrous}  There exists {$\ell \notin \bigcup_{i=0}^{\alpha-1} V_i \, \mod \, n$  if and only if}   $C(n,q\alpha)$ and $\gamma^{(\ell+1) \overline {q\beta}}\widetilde C(n,q\beta)$ are perfectly superimposable.
\end{corollary}

\Proof  By Proposition \ref{unionInter}, the union of the $V_i$'s is not a complete set of residues modulo $n$ if and only if $C(n,q\alpha)$ and $C(n,q\beta)$ are superimposable. Since $C(n,q\alpha)$ is the exact superimposition of the following $\alpha$ conjugates of $C(n,q)$
$$C(n,q), \gamma^{-\overline {q\alpha}}C(n,q),\ldots,\gamma^{-(\alpha-1)\overline{q\alpha}}C(n,q),$$
using the proof of Proposition \ref{unionInter} we get that  $\gamma^{-i\overline{q\alpha}}C(n,q)$ is perfectly superimposable with  $\gamma^{\ell \overline{q\beta}}C(n,q\beta)$ if and only if there exists a  $\ell \notin [-(q-1)\beta, q\beta[\, + i\overline \alpha \beta$ for {all} $0\leq i < \alpha$. Hence,  $C(n,q\alpha)$ and $\gamma^{\ell \overline{q\beta}}C(n,q\beta)$ are perfectly superimposable if and only if there exists $\ell \notin [-(q-1)\beta, q\beta[\, + i\overline \alpha \beta$ for {all} $0\leq i < \alpha$.  Finally, using Lemma \ref{lemMiroir},  $C(n,q\beta)=\gamma^{\overline {q\beta}}\widetilde C(n,q\beta)$ and we get that $C(n,q\alpha)$ and $\gamma^{\ell \overline{q\beta}}C(n,q\beta)$ are perfectly superimposable if and only if $C(n, q\alpha)$ and $\gamma^{\ell \overline{q\beta}}\gamma^{\overline {q\beta}}\widetilde C(n,q\beta)=\gamma^{(\ell +1)\overline{q\beta}}\widetilde C(n,q\beta)$ are so.~\QED

\begin{lemma} \label{propxyz} Let  $\alpha, \beta \in \N-\{0\}$, with $\alpha \perp \beta$, and let 
\begin{eqnarray}
x\alpha +y\beta&=&n-2\alpha \beta(q-1), \label{eqn3}
\end{eqnarray} with $q,\alpha, \beta \perp n$ and $q \geq 1$. Then:

\begin{itemize}
\item [ i)] Equation (\ref{eqn3}) always has a  solution $\{x, y\} \in \Z^2$;
\item [ ii)] it always has a unique solution with $1 \leq y \leq \alpha$;
\item [ iii)] if Equation (\ref{eqn3}) is satisfied, then $\alpha \perp (\alpha-y)$.
\end{itemize}
\end{lemma}

\Proof Since $\alpha \perp \beta$, {Bezout theorem yields that Equation} (\ref{eqn3}) always has a solution $\{x, y\} \in \Z^2$. Let us now suppose that there exist $2$ solutions, $\{x,y\}$ and $\{x',y'\}$, such that $1\leq y,y' \leq \alpha$.  Then $x\alpha +y\beta = x'\alpha +y' \beta$ and consequently, $\alpha(x-x')=\beta(y'-y)$. But $\alpha \perp \beta$ implies that  $\alpha | (y'-y)$: impossible, since $1 \leq y,y' \leq \alpha$. Finally, Equation (\ref{eqn3}) can be rewritten as
$$\alpha(x +2\beta(q-1))=n-y\beta$$
and since $\alpha \perp n$,  it follows that $\alpha \perp y$ and hence $\alpha \perp  (\alpha -y)$.  \QED

\begin{nota} \label{notayz} In the sequel, let $z=\alpha -y$, {where $y$ refers to the solution of Equation} (\ref{eqn3}). Let $i \in [0, \alpha-1]$ be one of the possible values of  $z$, as $z=\alpha-y$ and  $y \in [1,\alpha]$. Since $\alpha \perp z$ (see Lemma \ref{propxyz} iii)), {following Simpson} \cite{rjs2004}, there exists a unique $r(i) \in \N$ such that  $i\equiv  r(i)z \, \mod \, \alpha$. For $0\leq r < \alpha$, let
\begin{equation} \label{eqM}M(r)=r(x+(2q-1)\beta)-\left \lfloor \frac{zr}{\alpha}\right \rfloor \beta.\end{equation}
\end{nota}

The functions $r(i)$ and $M(r(i))$ will be useful in what follows, in order to obtain a new order for the  intervals $V_i$.

\begin{remark} Let $a=bq+r$, the Euclidean division of $a$ by $b$, with  $r<b$ and $a,b,q,r \in \N$. We have $r = a \, \mod \,  b$ and $\displaystyle q=\left \lfloor \frac{a}{b}\right \rfloor$. Thus, 
$$a=bq+r \iff a-bq-r=0 \iff a -b\left  \lfloor \frac{a}{b}\right \rfloor -(a\, \mod \, b) = 0.$$
\end{remark}

\begin{lemma} [Lemma 7 in \cite{rjs2004}] \label{lem7Simpson} For  $i \in [0,  \alpha-1]$, $M(r(i))\equiv -i\overline \alpha \beta \, \, \mod \, n.$
\end{lemma}

\Proof For a fixed $i$, let us consider  $\alpha(M(r(i))+i\overline \alpha \beta)$. In what follows, we will write $r$ instead of $r(i)$, in order to simplify the notation.  Using Equation (\ref{eqn3}) and the definition of $M(r)$, we get:
\begin{eqnarray*}
\alpha \left (M(r)+i\overline \alpha \beta\right )&\equiv &
\beta \left(zr-\alpha \left \lfloor \frac{zr}{\alpha} \right \rfloor-i\right) \, \mod \, n.
\end{eqnarray*}
{Replacing $a$ and $b$ in the previous remark by respectively  $zr$ and $\alpha$ and using the fact that $i \equiv zr \, \mod\,  \alpha$ yields that the term in parenthesis} has value $0$ and consequently, that $\alpha (M(r)+i\overline \alpha \beta) \equiv 0 \, \mod \, n$. Since $\alpha \perp n$,  $M(r)+i\overline \alpha \beta \equiv 0 \, \mod \,  n$ and we conclude. \QED

Lemmas \ref{lemInt}, \ref{lem311}, \ref{0ou1}, \ref{lemM} and \ref{lemfor} are not original results, since they appeared without emphasis in the proof of Theorem $8$ in \cite{rjs2004}. However, they are the key for the proofs of the results in the next section. 

\begin{lemma} \label{lemInt} Let $n \in \N$ be a fixed integer and let $I_0, I_1, \ldots, I_{r-1}$ be $r$ finite  intervals having same length and satisfying:
\begin{itemize}
\item [i)] $\max (I_0) - \min(I_{r-1})\geq n-1 \geq 1;$
\item [ii)] for  $0 \leq j <r-1$, if $I_{j+1}$ is located at the left of $I_j$, then $I_{j+1}\cup I_j$ is an interval.
\end{itemize}
Then $\bigcup_{j=0}^{r-1} I_j$ is a complete set of residues modulo $n$.
\end{lemma}

\Proof Let us suppose that the interval $I_{r-1}$ is not located at the left of the interval $I_0$. Since  $\max(I_0)-\min(I_{r-1}) \geq n-1$, it implies that $I_0 \cup I_{r-1}$ is an interval and that $I_0 \cap I_{r-1}=[\min(I_{r-1}), \max(I_0)]$. It follows that  $\Card(I_0 \cap I_{r-1})\geq n$ and that  $\bigcup_{j=0}^{r-1} I_j$ is a complete set of residues modulo $n$.
Let us now suppose  that the interval $I_{r-1}$ is located at the left of the interval $I_0$. By ii), there exist consecutive intervals that are located one to the left of the others. Condition ii) also insures  that all the integers between $I_{r-1}$ and $I_0$ are in the union of the $j$ intervals. Since $\max (I_0) - \min (I_{r-1}) \geq n-1$, the number of integers between the beginning of the interval $I_{r-1}$ and the end of the interval $I_0$ is at least $n$. In both cases, $\bigcup _{j=0}^{r-1} I_j$ is a complete set of residues modulo $n$. \QED

\begin{lemma} \label{lem311} Let $I_0, I_1, \ldots, I_{r-1}$ be finite intervals in $\Z$ and let $I$ be the shortest interval that contains {them}. Let us suppose that
\begin{itemize}
\item [i)] $I\setminus \bigcup _{j=0}^{r-1} I_j$  is non-empty;
\item [ii)] if $x \in \bigcup_{j=0}^{r-1} I_j$ and $y \in I \setminus \bigcup_{j=0}^{r-1} I_j$, then $|y-x|<n$.
\end{itemize}
Then $\bigcup_{j=0}^{r-1}I_j$ does not contain all the integers modulo $n$. 
\end{lemma}

\Proof Let us suppose that $y \equiv x \, \mod \, n$, for $x, y$ satisfying condition ii). By condition i), such a $y$ exists. Then there also exists $k \in \Z$ such that $y-x =kn$. Since $y \notin \bigcup_{j=0}^{r-1} I_j$ and $x \in \bigcup_{j=0}^{r-1} I_j$, we have $k\neq 0$, otherwise $x=y$. Thus $|k| \geq 1 \Rightarrow |y-x| \geq n$, which contradicts ii).  \QED

\begin{lemma} \label{0ou1}  Let $z=\alpha-y$ and $0 \leq r < \alpha$ {be such as in Notation} \ref{notayz}.  Then
$\left \lfloor\frac{ z(r+1)}{\alpha} \right \rfloor - \left \lfloor \frac{zr}{\alpha }\right \rfloor \in \{0,1\}.$ 
\end{lemma}

\Proof Let $zr=i\alpha +t$, with $i \in \N$ and $0 \leq t < \alpha$. Then 
\begin{eqnarray*}
\left \lfloor\frac{ z(r+1)}{\alpha}  \right \rfloor - \left \lfloor \frac{zr}{\alpha }\right \rfloor &=& 
\left \lfloor\frac{ i\alpha +t +z}{\alpha} \right \rfloor - \left \lfloor \frac{i\alpha +t}{\alpha }\right \rfloor  = i+ \left \lfloor \frac{t +z}{\alpha} \right \rfloor - i- \left \lfloor \frac{t}{\alpha }\right \rfloor = \left \lfloor \frac{t +z}{\alpha} \right \rfloor - \left \lfloor \frac{t}{\alpha }\right \rfloor . 
\end{eqnarray*}
Since $1 \leq y \leq \alpha$ and $0 \leq t < \alpha$, we have $0\leq t+z < 2\alpha$ and consequently, 
$$ \left \lfloor \frac{t +z}{\alpha} \right \rfloor - \left \lfloor \frac{t}{\alpha }\right \rfloor = \left \lfloor \frac{t +z}{\alpha} \right \rfloor - 0 \leq 1.$$ \vspace{-1cm}{\flushright \QED}

\begin{lemma} \label{lemM} Let $M(r)$ be defined as in Equation (\ref{eqM}). Then 
\begin{itemize}
\item [i)] $M(0)=0$;
\item [ii)] $M(\alpha -1)=n-x-2\beta(q-1)-i$, with $i=0$ if $y \neq \alpha$ and $i=\beta$ otherwise.
\end{itemize}
\end{lemma}

\Proof 
\begin{itemize}
\item [i)] $M(0)=0(x+(2q-1)\beta)-\left \lfloor \frac{z\cdot 0}{\alpha}\right \rfloor \beta=0.$
\item [ii)] If $y \neq \alpha$, then
\begin{eqnarray}
M(\alpha-1)&=&(\alpha-1)(x+(2q-1)\beta)-\left \lfloor \frac{z(\alpha-1)}{\alpha} \right \rfloor \beta \label{16} \\ 
& =& \alpha x+(2q-1)\alpha \beta-x-(2q-1)\beta-z\beta+\beta \label{17}\\
&=&\alpha x+y\beta -\alpha\beta+(2q-1)\alpha \beta-x+\beta-(2q-1)\beta \nonumber\\
&=& n-x-2\beta(q-1). \label{18}
\end{eqnarray}
Equation (\ref{17}) is deduced from Equation (\ref{16})  using the fact that
$$\left \lfloor \frac{z(\alpha-1)}{\alpha}\right \rfloor \beta = z\beta +\left \lfloor \frac{-z}{\alpha} \right \rfloor \beta = z\beta - \beta ,$$ 
since $0 < z=\alpha -y < \alpha$, as $1\leq y < \alpha$.\\
 If $y=\alpha$, then $z=0$ and 
$$\left \lfloor \frac{z(\alpha-1)}{\alpha}\right \rfloor \beta= 0,$$ 
implying 
\begin{equation}
M(\alpha-1)= n-x-2\beta(q-1)-\beta. \label{19}
\end{equation}  \vspace{-1cm}\flushright{\QED} 
\end{itemize}

\begin{lemma} \label{lemfor} Let $M(r)$ be defined as in Equation (\ref{eqM}). Then
\begin{itemize}
\item [i)] $M(r+1)-M(r)= x+(2q-1)\beta -\beta \left ( \left \lfloor \frac{z(r+1)}{\alpha}\right \rfloor - \left \lfloor \frac{zr}{\alpha} \right \rfloor \right )$;
\item [ii)] if $x \leq 0$, then $M(r+1)-M(r) \leq \beta(2q-1)$.
\end{itemize}
\end{lemma}

\Proof We have:
\begin{eqnarray*}
M(r+1)-M(r)&=& \left ((r+1)(x+(2q-1)\beta)-\left \lfloor \frac{z(r+1)}{\alpha}\right \rfloor \beta \right ) -\left (r(x+(2q-1)\beta)-\left \lfloor \frac{zr}{\alpha}\right \rfloor \beta \right) \nonumber \\
&=& x+(2q-1)\beta -\beta \left ( \left \lfloor \frac{z(r+1)}{\alpha}\right \rfloor - \left \lfloor \frac{zr}{\alpha} \right \rfloor \right ),
\end{eqnarray*}
which is $ \leq \beta(2q-1) $ if $x \leq 0$, using Lemma \ref{0ou1}.  \QED

{The following theorem is a particular case of Theorem 8 in} \cite{rjs2004} {which first appeared in} \cite{rm1985}. 
\begin{theorem} [\cite{rm1985,rjs2004}]\label{thprincipal} $C(n,q\alpha)$ and $C(n,q\beta)$ are superimposable if and only if there {exists}  $\{x, y\}\in \{\N-\{0\}\}^2$ such that  
\begin{eqnarray}
x\alpha + y \beta  &=& n-2\alpha \beta(q-1). \label{eqn4}
\end{eqnarray}
\end{theorem}

\Proof For the Christoffel words $C(n,q\alpha)$ and $C(n,q\beta)$, by Lemma \ref{propxyz}  there {exists} a unique $\{x, y\}$ satisfying {Equation} (\ref{eqn4}), with $1\leq y \leq \alpha$.  We want to show that  $C(n,q\alpha)$ and $C(n,q\beta)$ are superimposable if and only if $x > 0$. 

($\Longrightarrow$) Let us suppose that $x\leq0$ and let us consider the union of the intervals given in Equation (\ref{unionVi}). Using Lemma \ref{lem7Simpson}, we have, modulo $n$,
\begin{eqnarray} \label{eqlala} 
 \bigcup _{i=0}^{\alpha-1} \left  \{[(-q+1)\beta, q\beta -1] +i\overline \alpha \beta \right \}&=&\displaystyle \bigcup _{r=0}^{\alpha-1} \left  \{[(-q+1)\beta, q\beta -1] -M(r) \right \}.
\end{eqnarray}
Let $I_r=[(-q+1)\beta, q\beta-1]-M(r)$, for $0 \leq r < \alpha$.  Then using Lemma \ref{lemM}, we get $\max(I_0)=q\beta -1 -M(0)=q\beta-1$ and $\min(I_{\alpha-1})=(-q+1)\beta -M(\alpha-1)$ and then
\begin{eqnarray*}
\max(I_0)-\min(I_{\alpha-1})&=&q\beta -1- \left((-q+1)\beta -M(\alpha-1)\right )\\
&=& q\beta -1 +q\beta -\beta +n-x-2\beta(q-1)-i\\
&=& \beta +n-x-1-i,
\end{eqnarray*}
where $i=0$ if $y\neq \alpha$ and $i=\beta$ otherwise (see Lemma \ref{lemM}).
Since $x \leq 0$, $-x$ is  non-negative. Hence $\max(I_0)-\min(I_{\alpha-1}) \geq n-1$. 

If  $I_r \cup I_{r+1}$ is not an interval then  $M(r+1)-M(r)>|I_r|+1$. Thus, in order to show that $I_r\cup I_{r+1}$ is an interval, it is sufficient to show that $M(r+1)-M(r) \leq |I_r|+1$. By Lemma \ref{lemfor} ii), we have $M(r+1)-M(r) \leq \beta(2q-1)$. Moreover, all the intervals have length 
$$|I_r|=q\beta-1-(-q+1)\beta=2q\beta-\beta-1=\beta(2q-1)-1.$$ Hence $M(r+1)-M(r)\leq |I_r|+1$. Then, for all $0 \leq r < \alpha-1 $, $I_r \cup I_{r+1}$ is an interval. 
Recall that Proposition \ref{unionInter} tells us that $C(n,q\alpha)$ and $C(n, q\beta)$ are superimposable if and only if the union (\ref{eqlala}) is not a complete set of residues modulo $n$. Applying  Lemma \ref{lemInt}, we conclude that for  $x\leq 0$  the words are not superimposable.

($\Longleftarrow$) Let us now suppose that $x>0$ and let us show that it implies that the words are superimposable. By Proposition \ref{unionInter}, it is sufficient to show that if  $x>0$, then $\bigcup _{r=0}^{\alpha -1}\left \{ [-(q-1)\beta, q\beta -1] -M(r)\right \}$ does not contain all the integers modulo $n$. 

Let us recall that $I_r=[-(q-1)\beta,q\beta-1]-M(r)$, for $0 \leq r < \alpha$. Since $x>0$ and $q\geq 1$, and using Lemmas \ref{0ou1} and \ref{lemfor}, we have
\begin{eqnarray*}
M(r+1)-M(r)&=& x+(2q-1)\beta-\beta \left (\left \lfloor\frac{ z(r+1)}{\alpha} \right \rfloor - \left \lfloor \frac{zr}{\alpha }\right \rfloor \right)\\
&\geq&x+(2q-1)\beta -\beta\\
&=& x+2\beta(q-1) \geq x >0.
\end{eqnarray*}

Thus, the intervals $I_r$ are located one to the left of the others, for $0 \leq r< \alpha$. They all have the same cardinality, that is: $\Card(I_r) =|I_r| +1 = \beta(2q-1)-1+1= \beta(2q-1).$

Let us suppose that  $I=\bigcup_{r=0}^{\alpha-1} I_r$ is not an interval. Then condition i) of Lemma \ref{lem311} is satisfied. For condition ii), it is sufficient to take $y= \min (I_0) -1$ (since $\max(I \setminus \cup I_j) \leq \min (I_0)-1$) and $x=\min(I)$ and to check that $y-x < n$. We have $x= -(q-1)\beta-M(\alpha-1)$ and $y=-(q-1)\beta-1$.

Consequently: 
$$y-x=\left (-(q-1)\beta-1\right )-\left(-(q-1)\beta-\left(n-x-2\beta(q-1)-i\right)\right ) = n-x-1-2\beta(q-1)-i,$$ with $i\in \{0,\beta\}$. 
Since $x>0$, this value is $<n$. By Lemma \ref{lem311}, we conclude that the union of these intervals does not contain all the integers modulo $n$. \QED

\subsection{Number of superimpositions of Christoffel words} \label{secnbsuperp}

In this section, we prove the exact number of superimpositions of two Christoffel words having same length and we give a shift that always allows a perfect superimposition for two superimposable Christoffel words.

\begin{defn} Let $C(n,\alpha)$ and $C(n,\beta)$ be two superimposable Christoffel words. The  {\it number of superimpositions} of these two words is defined by 
$$\Card \left (\{k \in [0,n-1] \,  | \, C(n,\alpha) \textnormal{ and } \gamma^{k}C(n,\beta) \textnormal{ are exactly superimposable}\}\right ).$$
\end{defn}

Some results are first required.

\begin{corollary} [of Lemma \ref{lem311} and of its proof] \label{desctrous}If the two conditions of  Lemma \ref{lem311} are satisfied, then
\begin{itemize}
\item [i)] the elements of $ I \setminus \bigcup_{j=0}^{r-1} I_j$ are all distinct modulo $n$;
\item [ii)] if $\Card(I) \geq n$, then modulo $n$, the elements of $I \setminus \bigcup_{j=0}^{r-1} I_j$ are exactly the ones that are not in $\bigcup_{j=0}^{r-1} I_j$;
\item [iii)] if $\Card(I) < n$, then modulo $n$, the elements of $\Z \setminus \bigcup_{j=0}^{r-1} I_j$ are exactly the ones that are not in $ \bigcup_{j=0}^{r-1} I_j$ and the following $n-\Card(I)$ elements:
$$\{ \min(I)-(n-\Card(I)), \ldots, \min(I)-2, \min(I)-1\}.$$
\end{itemize}
\end{corollary}

\Proof \begin{itemize}
\item [{i)}] Let $x,y \in I \setminus \bigcup_{j=0}^{r-1} I_j$ and without loss of generality, let us suppose that $y>x$. Then $y \leq \max (I \setminus \bigcup_{j=0}^{r-1} I_j)$ and $x > \min(I)$, since $\min(I)$ is contained in  $\bigcup_{j=0}^{r-1} I_j$. Consequently $y-x < \max(I \setminus \bigcup_{j=0}^{r-1} I_j)-\min(I)$ which is, by Lemma \ref{lem311} ii), $< n$. Hence $y-x< n$.

\item [{ii)}] Since $I$ is an interval and $\Card(I)\geq n$, $I$ contains all the elements $\mod \, n$. By Lemma \ref{lem311} ii), there is no element  in $\bigcup_{j=0}^{r-1} I_j$ that is equal, modulo $n$, to an element in $I \setminus \bigcup_{j=0}^{r-1} I_j$. 

\item [{iii)}] One can easily observe that the $n-\Card(I)$ elements are not equal, modulo $n$, to any element of $I$. \QED
\end{itemize}

\begin{lemma} \label{lemPartieEntiere} Let $C(n, j) \in \{a<b\}^*$ be a Christoffel word. Then 
$${\displaystyle C(n,j)[i]= \left \{ \begin{array}{ll} 
a & \textnormal{if } \displaystyle \left \lfloor \frac{n-j}{n} (i+1)\right \rfloor - \left \lfloor \frac{n-j}{n} i\right \rfloor =0 \\ \\
b & \textnormal{if }  \displaystyle \left \lfloor \frac{n-j}{n} (i+1)\right \rfloor - \left \lfloor \frac{n-j}{n} i\right \rfloor =1.
\end{array}
      \right . }$$
\end{lemma}

\Proof Follows from the definition of Christoffel words. Doing the difference between the integer parts corresponds to check if there is a multiple of $n$ or not between both values. If the difference is $0$, then no multiple of $n$ occurs. \QED

\begin{proposition} \label{nbsup} Let $C(n,q\alpha)$ and $C(n, q\beta)$ be two superimposable Christoffel words. The number of superimpositions of $C(n,q\alpha)$ and $C(n,q\beta)$ is
\begin{itemize}
\item [i)] $xy$, if $x \leq \beta$;
\item [ii)] $x\alpha+y\beta-\alpha\beta$, if $x>  \beta$;
\end{itemize}
where $\{x, y\}$ is the unique solution of Equation (\ref{eqn4}), with $1\leq y \leq \alpha$. 
\end{proposition}

\Proof Let us recall from Theorem \ref{thprincipal} that if two Christoffel words are superimposable and if the solution of Equation (\ref{eqn4}) is $\{x,y\}$ with $1\leq y \leq \alpha$, then $x>0$. Let us denote by $I$ the shortest interval that contains the union of the intervals given in Equation (\ref{eqlala}). Then using Lemma \ref{lemM}, we get
\begin{eqnarray*}
\Card(I)&=& \max(I)-\min(I)+1 \\
&=& \max(I_0)-\min(I_{\alpha-1}) +1\\
&=& (q\beta-1)-((-q+1)\beta-M(\alpha-1)) +1 \\
&=& q\beta -1 +q\beta - \beta +(n-x-2\beta(q-1)-i)+1\\
&=& n-x+\beta-i,
\end{eqnarray*}
with $i=0$ if $y \neq \alpha$, and $i=\beta$ otherwise. 

\begin{itemize}
\item [i)] \underline {Let us suppose that $x\leq \beta$ and $y\neq \alpha$.} Then $\Card(I) = n-x+\beta \geq n$. By  Corollary \ref{desctrous} ii), the complementary set modulo $n$ of $\bigcup_{j=0}^{\alpha-1} I_j$ has the same  cardinality as the number of elements contained between $I_0$ and $I_1$, $I_1$ and $I_2$, etc. The number of elements contained between $I_r$ and $I_{r+1}$ is $M(r+1)-M(r)-(2q-1)\beta$, that is the distance between the beginning of both intervals minus the cardinality of one interval. Using Lemma \ref{lemfor}, we have
\begin{eqnarray}
M(r+1)-M(r)-(2q-1)\beta&=&x-\beta\left (\left \lfloor \frac{z(r+1)}{\alpha} \right \rfloor -\left \lfloor\frac{ zr}{\alpha} \right \rfloor \right ). \label{nbelts}
\end{eqnarray}
There is a gap between two intervals if the value of (\ref{nbelts}) is $>0$. This value corresponds to the number of integers contained in the gap.
Since $x \leq \beta$, it will be the case for all $r$ such that
$ \left \lfloor \frac{z(r+1)}{\alpha} \right \rfloor -\left \lfloor\frac{ zr}{\alpha} \right \rfloor=0$.
Using  Lemma \ref{lemPartieEntiere}, with $j=y$, $i=r$ and $n=\alpha$, we find that it is the case for exactly $y$ values of $r$. Thus, there are $xy$ possible superimpositions.

\item [ii)] \underline {Let us suppose that $x > \beta$.} One can easily observe that $x>\beta \Longrightarrow y \neq \alpha$. Then $\Card(I) = n-x+\beta < n$. We still have that  
$ \left \lfloor \frac{z(r+1)}{\alpha} \right \rfloor -\left \lfloor\frac{ zr}{\alpha} \right \rfloor=0$ for $y$ values of  $r$. Moreover, since $0 \leq r < \alpha$, there are $(\alpha-1)$ gaps containing each $$x-\beta\left (\left \lfloor \frac{z(r+1)}{\alpha} \right \rfloor -\left \lfloor\frac{ zr}{\alpha} \right \rfloor \right )$$ integers. Thus, by Lemma \ref{lemPartieEntiere}, 
$$\left \lfloor \frac{z(r+1)}{\alpha} \right \rfloor -\left \lfloor\frac{ zr}{\alpha} \right \rfloor=1$$ for $\alpha-1-y=z-1$ values. Hence $x-\beta\left (\left \lfloor \frac{z(r+1)}{\alpha} \right \rfloor -\left \lfloor\frac{ zr}{\alpha} \right \rfloor \right )=x-\beta$ for $(z-1)$ values of $r$. Using Corollary \ref{desctrous} iii), we know that there are $n-\Card(I)$ others possible values outside the interval $I$. The number of superimpositions is then given by
\begin{eqnarray*}
xy+(x-\beta)(z-1)+ n-\Card(I) &=& xy+(x-\beta)(\alpha-y-1)+n-(n-x+\beta)\\
&=& x\alpha +y\beta -\alpha \beta.
\end{eqnarray*}

\item [iii)] \underline {Let us suppose that $x\leq \beta$ and $y=\alpha$.} Then $\Card(I)= n-x < n$. This case is similar to case ii), except that here, since $y=\alpha$, 
$$ \left \lfloor \frac{z(r+1)}{\alpha} \right \rfloor -\left \lfloor\frac{ zr}{\alpha} \right \rfloor=0$$ between every interval and hence $z=\alpha-y=0$. Since there are $(y-1)=(\alpha-1)$ gaps between $I_0$ and $I_{\alpha -1}$, using Corollary \ref{desctrous} iii), we find that the number of possible superimpositions is given by 
$$x(y-1)+n-\Card(I)=x(y-1)+n-(n-x)=xy.$$
\vspace{-1.8cm}{\flushright \QED}
\end{itemize}

\begin{remark} Since Lemma \ref{propxyz}, we have supposed that  $\{x,y\}$  is the solution of Equation  (\ref{eqn4}) such that $y\leq \alpha$. In the proof of Proposition \ref{nbsup}, we still use this assumption. It is possible to rewrite all these results considering the solution for which $x\leq \beta$. We would have obtained similar result as in Proposition \ref{nbsup}, with the conditions  $y \leq \alpha$ and $y>\alpha$. 
\end{remark}

Theorem \ref{uneSuperpGen} is a generalization of  Corollary \ref{uneSuperp} for any values of  $q, \alpha, \beta$, such that $\alpha \perp \beta$: for two superimposable Christoffel words having same length, we give a shift that always allows a perfect superimposition.

\begin{theorem} \label{uneSuperpGen} Let $C(n,q\alpha)$ and $C(n,q\beta)$ be two superimposable Christoffel words, with $\alpha \perp \beta$. Then $C(n,q\alpha)$ and $\gamma^{1-r}\widetilde C(n,q\beta)$ are perfectly superimposable, where $qr\equiv 1 \, \mod \, n$.
\end{theorem}

\Proof By Corollary \ref{corTrous}, $C(n,q\alpha)$ and $\gamma^{(\ell+1)\overline{q\beta}}\widetilde C(n,q\beta)$ are superimposable if and only if $\exists \ell \notin  \bigcup _{i=0}^{\alpha-1}V_i \, \mod \, n$. It is then sufficient to show that there exists $\ell\notin \bigcup _{i=0}^{\alpha-1}V_i \, \mod \, n$ such that $(\ell +1)\overline {q \beta} \equiv 1-r$. Isolating $\ell$, we get that this last condition is equivalent to 
$$\ell \equiv -q\beta +rq\beta -1\equiv \beta-1-q\beta.$$

Let us show that $\beta -1 -q\beta \notin \bigcup_{i=0}^{\alpha-1} V_i \, \mod \, n$.  

If $\alpha=1$, by Equation (\ref{eqlala}), the union of the $V_i$'s is the interval $[-q\beta+\beta, q\beta-1]$. Then $\beta-q\beta-1$ is the element preceding the interval and since the words are superimposable, the interval has a length  $< n$, and consequently $\beta - q\beta -1 \, \mod \, n$ is not contained in the interval.  

If $\alpha >1$, let us consider the intervals $I_0$ and $I_1$. There exist elements between both intervals, since 
$$\displaystyle M(1)-M(0)-(2q-1)\beta=x+(2q-1)\beta - \left \lfloor \frac{z}{\alpha}\right \rfloor\beta-0-(2q-1)\beta=x>0,$$ 
as $z=\alpha-y < \alpha$. Moreover, 
\begin{eqnarray}
]\max(I_1), \min(I_0)[ &=& \left ]q\beta -1 -\left (x+(2q-1)\beta\right), (-q+1)\beta\right[\\
& =&\left ]q\beta -1 -x-2q\beta +\beta, -q\beta+\beta\right[ \nonumber \\
&=& ]\beta -q\beta -1-x, -q\beta +\beta[. \label{trou}
\end{eqnarray}
Thus, $\beta-q\beta-1$ is located between $I_0$ and $I_1$. In order to conclude, it is sufficient to show that this element does not appear in an other interval. It is true if $(\beta-q\beta-1)-\min(I_{\alpha-1})<n$. Let us verify:
\begin{eqnarray*}
(\beta-q\beta-1)-\min(I_{\alpha-1}) &=& \beta-q\beta-1-((-q+1)\beta-(n-x-2\beta(q-1)-i))\\
&=& n-2\beta(q-1) -x-1-i < n 
\end{eqnarray*}
where $i\in \{0,\beta\}$.  \QED

\section{Generalization to words having different lengths}

In this section, we use Theorem \ref{th2Simpson} in order to generalize the results of Sections \ref{secmemelng} and \ref{secnbsuperp} for arbitrary Christoffel words, not necessarily having same length.

\begin{theorem} {\rm \cite{rm1985, rjs2004}}\label{theoremprinc} Let  $C(n,q\alpha)$ and $C(m,q\beta)$ be Christoffel words, with $\alpha \perp \beta$. Then $C(n,q\alpha)$ and $C(m,q\beta)$ are superimposable if and only if there {exists}  $\{x,y\} \in \{\N-\{0\}\}^2$ such that
\begin{eqnarray}
x\alpha +y\beta &=&p -2\alpha\beta(q-1), \label{eqnsatis}
\end{eqnarray}
with $p=\gcd(m,n)$. 
\end{theorem}

\Proof By Theorem \ref{th2Simpson},  $C(n,q\alpha)$ and $C(m,q\beta)$ are superimposable if and only if $C(p,q\alpha)$ and $C(p,q\beta)$ are so. We conclude using Theorem \ref{thprincipal}, since it insures that $C(p,q\alpha)$ and $C(p,q\beta)$ are superimposable if and only if there {exists} $\{x, y\} \in \{\N-\{0\}\}^2$ satisfying Equation (\ref{eqnsatis}). \QED

\begin{lemma} \label{remsuperp} If $C(p,\alpha)$ and $\gamma^{-k}C(p, \beta)$ are perfectly superimposable, then $C(n,\alpha)$ and $\gamma^{-k+ip} C(m,\beta)$ are so, with $m>n$,  $p=\gcd(n,m)$ and $0\leq i< \frac{m}{p}$. 
\end{lemma}

\Proof Theorem \ref{th2Simpson} shows that $C(p,\alpha)$ and $\gamma^{-k}C(p, \beta)$ are perfectly superimposable if and only if $C(n,\alpha)$ and $\gamma^{-k} C(m,\beta)$ are so. Moreover, one can easily observe that  $C(p,\alpha)$ and $\gamma^{-k}C(p, \beta)$ are perfectly superimposable if and only if $C(p,\alpha)$ and $\gamma^{-k+ip} C(p,\beta)$ are so. These $-k+ip$ correspond to different shifts, for $0\leq i< \frac{m}{p}$, for words of length at most $m$. \QED

\begin{proposition} Let $C(n,q\alpha)$ and $C(m,q\beta)$ be two superimposable Christoffel words, with $\alpha \perp \beta$, $p=\gcd(m,n)$ and $m>n$. The number of superimpositions is
\begin{itemize}
\item [i)]  $\displaystyle xy\frac{m}{p}$, if $x\leq \beta$;
\item [ii)] $\displaystyle (x\alpha+y\beta-\alpha\beta)\frac{m}{p}$, if $x > \beta$;
\end{itemize}
with $\{x, y\}\in \{\N-\{0\}\}^2$ the solution of $x\alpha+y\beta=p-\alpha \beta (q-1)$ such that  $y \leq \alpha$.
\end{proposition}

\Proof Follows from Proposition \ref{nbsup} and from Lemma \ref{remsuperp}. \QED

\begin{theorem} Let  $C(n,q\alpha)$ and $C(m,q\beta)$ be two superimposable Christoffel words, with $\alpha \perp \beta$ and $p=\gcd(m,n)$. Then $C(n,q\alpha)$ and $\gamma^{-(r-1)+ip}\widetilde C(m,q\beta)$ are perfectly superimposable, with $qr\equiv  1 \,\mod \, p$ and $\displaystyle 0 \leq i < \frac{m}{p}$. 
\end{theorem}

\Proof Follows from Theorem \ref{uneSuperpGen} and from Lemma \ref{remsuperp}. \QED

\section{Other results}

In this last section, we first give a new necessary and sufficient condition for the perfect superimposition of two Christoffel words $C(n,\alpha)$ and $C(n,\beta)$, with $\alpha \perp \beta$. Then we give a result concerning the word obtained by the superimposition of two Christoffel words having same length. We end this section by a new proof of a problem related to the money problem, using the geometric interpretation of Christoffel words. 

\begin{theorem} \label{theore51}Let $u=C(n,\alpha) \in \{a< z\}^n$ and $v =C(n,\beta)\in\{b< z\}^n$ be Christoffel words. There {exists} $\{x, y\} \in \{\N-\{0\}\}^2$ such that $\alpha x+\beta y=n$ if and only if $u$ and $\widetilde{v}$ are perfectly superimposable.
\end{theorem}

\Proof  
($\Longrightarrow$) Let us suppose that there {exists} $\{x, y\} \in \{\N-\{0\}\}^2$ such that $\alpha x + \beta y=n$. Let us consider the Christoffel words $u'=C(n,x\alpha)$ and $v'=C(n,y\beta)$. Since $\alpha x +\beta y=n$, these words are complementary, that means that $u'$ and $\widetilde v'$ are perfectly superimposable. Using Lemma \ref{lempos}, we conclude that $u$ and $\widetilde v$ are so.

($\Longleftarrow$) Let us suppose that $u$ and $\widetilde{v}$ are perfectly superimposable. Let $d=\gcd(\alpha, \beta)$. By Lemma \ref{lempos}, $C(n,d)\in \{a< z\}^n$ and $\widetilde C(n,d) \in \{z<b\}^n$ are also superimposable. Let us now show that $u$ and $\widetilde{v}$ are superimposable only if  $d | n$. If $d\not |n$, then $C(n,d)$ can be written as the product of  $az^i$ and $az^{i+1}$, it begins by  $az^i$ and ends by $az^{i+1}$. Moreover, $\widetilde{C}(n,d)$ ends by $bz^ib$. There is a conflict between a letter $a$ and a letter $b$, since
$$C(n,d)=paz^iz \, \, \textnormal{and} \, \, \widetilde{C}(n,d)=p'bz^i b.$$ 
Thus, if the words are perfectly superimposable,  $d | n$. 

Moreover, since $d=\gcd(\alpha, \beta)$, $d|\alpha$ and $d|\beta$. Thus, $\gcd(n,\alpha)=d$ and $\gcd(n,\beta)=d$. Since $C(n,\alpha)$ and $C(n,\beta)$ are Christoffel words, $\alpha \perp n$ and $\beta \perp n$. Hence $d=1$. Applying Theorem \ref{theoremprinc} with $m=n$, $q=1$, we get  $p=1$ and consequently, $C(n,\alpha)$ and $C(n,\beta)$ are superimposable if and only if there {exists}  $\{x,y\} \in \{\N-\{0\}\}^2$ such that $x\alpha +y\beta=n$. \QED

\begin{defn} Let $w \in \{a,b\}^*$ and let $A=\{i_1, i_2, \ldots, i_k\}$ be the set of positions of the $a$'s in $w$. Then the word $w'\in \{a,b\}^*$ obtained by the  {\it decimation} $D_{p/q, a}$ of $w$ over the letter $a$, with $p\leq q$, is the word $w$ for which we have deleted the letters $w[i_j]$,  for all $j \in \{\ell q +1, \ell q +2, \ldots, \ell q +p\}_{0 \leq \ell \leq \lfloor |A|/q\rfloor}$ if $p/q <0$ and  for all $j \in \{|A|-\ell q, |A|-\ell q -1, \ldots, |A|-\ell q -p +1\}_{0 \leq \ell \leq \lfloor |A|/q\rfloor}$ otherwise. In other words, $w'$ is the word $w$ for which $p$ occurrences over $q$ of the letter $a$ are removed from left to right if $p/q<0$, and from right to left otherwise.
\end{defn}

{\ex Let consider $w=aabaabababa$. The decimation $D_{1/3, a}(w)$ yields $w'= abababab$. Then performing $D_{-1/2,b}$ over $w'$ gives $w''=aabaab$.
}

\begin{theorem}\label{theodec} Let $u=C(n,\alpha)\in \{a< z\}^n$ and $C(n,\beta)\in \{b<z\}^n$ be two superimposable Christoffel words with $\alpha \perp \beta$. Let $v$ be the conjugate of $C(n,\beta)$ that is perfectly superimposable to $u$. Let $w$ be defined as
$$ w[i]= \left \{ \begin{array}{ll} a & \textnormal{if }  u[i]=a\\
b & \textnormal{if }  v[i]=b \\
z & \textnormal{otherwise.}
\end{array}
      \right . $$
Let $w'$ be the word obtained from $w$, after having removed the letter $z$. Then $w'$ is the Christoffel word of slope  $\beta / \alpha$. 
\end{theorem}

\Proof Let $\{x, y\}\in \{\N-\{0\}\}^2$ be such that  $\alpha x+ \beta y=n$. By Theorem \ref{thprincipal}, we know that such $x,y$ exist. Let us consider the Christoffel word $t\in\{a<b\}^n$ with $\alpha x$ occurrences of the letter $a$ and $\beta y$ occurrences of the letter $b$. Let us perform the decimation $D_{(x-1)/x, a}(t)$: it removes $(\alpha x-\alpha)$ letters $a$'s. The decimation $D_{-(y-1)/y,b}$ over the word obtained removes  $(\beta y -\beta)$ letters $b$'s. Since the decimation operation preserves Christoffel words \cite{jpb2001}, the word $w'$ obtained is a Christoffel word of length $\alpha + \beta$ with $\alpha$ occurrences of the letter $a$ and $\beta$ occurrences of the letter $b$.  \QED

{\ex \label{exa22} {Let $u=C(13,4)=azzazzazzazzz$ and $C(13,3)=bzzzbzzzbzzzz$. These words are   superimposable. Indeed, it is sufficient to take the conjugate $v=\widetilde C(13,3)=zzzzbzzzbzzzb$. We then find $w=azzabzazbazzb$ and $w'=aababab$. Note that the equation $4x+3y=13$ has the solution $x=1$ and $y=3$. Thus, we consider the Christoffel word $t=C(13,4)=abbabbabbabbb$. The decimation $D_{{0/1},a}(t)$ does not erase any $a$. Then we perform $D_{-2/3, b}$ over $D_{{0/1},a}(t)$:  starting from the left we erase $2$ occurrences over $3$ of $b$'s. We get $aababab=C(7,4)$.}


\subsection{Money problem}

In {Theorem} \ref{theore51}, we showed that two Christoffel words $u$ and $\tilde v$ of length  $n$ are perfectly superimposable if and only if there exist integers  $\alpha, \beta$ such that $\alpha x+\beta y=n$. {In what follows,} $\alpha x + \beta y$ occurs again: we prove, using the geometric interpretation of Christoffel words, {classical} results of Sylvester concerning the money problem, also known as  {\it Frobenius problem}. 

Let us first recall the money problem.

\begin{defn} \cite{eww2007} Let $0 < a_1  < \ldots < a_n $ be $n$ integers, with $n\geq 2$, that represent  $n$ different values of money pieces and such that $\gcd(a_1,a_2,\ldots, a_n)=1$. The possible amounts of money that can be obtained using these $n$ pieces are given by
$$\sum_{i=1}^{n}a_ix_i,$$
where $x_i\in \N$ denotes the number of the piece $a_i$ used. The {\it money changing problem} consists of determine the greatest integer $N=g(a_1,a_2,\ldots, a_n)$ that cannot be obtained using the pieces of money  $a_1, a_2, \ldots, a_n$. This integer is called the  {\it  Frobenius number}. 
\end{defn}

If $a_1=1$, all amounts can be obtained. It is not the case in general: only a few amounts can be obtained. For instance, with pieces of  $2, 5$ and $10$, it is impossible to obtain $1$ and $3$, while all the other quantities can be obtained. Hence $g(2,5,10)=3$. 

\begin{proposition} \label{propSyl}{\rm \cite{jjs1884}} The greatest integer that cannot be obtained with the pieces $a$ and $b$ is 
\begin{eqnarray}
g(a,b)&=&(a-1)(b-1)-1.\label{sylvest}
\end{eqnarray}
\end{proposition}

Proposition \ref{prop2Syl} appears in \cite{eww2007}, but the origin is unknown.

\begin{proposition} \label{prop2Syl} The number of integers that cannot be obtained with the pieces $a$ and $b$ is given by
\begin{eqnarray}
\frac{(a-1)(b-1)}{2}. \label{nbsylvester}
\end{eqnarray}
\end{proposition}

\begin{corollary} \label{corMonnaie}The number of elements of the submonoid {of $N$} generated by $a$ and $b$ and smaller than  $(a-1)(b-1)$ is $\frac{(a-1)(b-1)}{2}.$
\end{corollary}

\Proof  We know by Proposition \ref{propSyl} that all the integers greater or equal to  $(a-1)(b-1)$ are representable with $a$ and $b$. Thus the unrepresentable $\frac{(a-1)(b-1)}{2}$ integers  given in  Proposition \ref{prop2Syl} are necessarily smaller than $(a-1)(b-1)$. Since half of the $(a-1)(b-1)$ elements smaller than $(a-1)(b-1)$ (including the 0) are not representable with $a$ and $b$, there is exactly the same quantity that is representable. \QED

In what follows, we will show that it is possible to prove Corollary \ref{corMonnaie} using the geometric representation of Christoffel words and their Cayley graphs. 

\begin{theorem} \label{thmbande} Let $a, b \in \N$. Let us consider the quadrant defined by $x \geq 0$ and $y \leq 0$, having at the coordinate $(x,-y)$ the value $xb+ya$. While considering only the integer coordinates $(x,-y)$ such that $xb+ya < ab$, the boundary obtained can be coded by a Christoffel word having exactly $a$ occurrences of the letter $\alpha$ and $b$ occurrences of the letter $\beta$. 
\end{theorem}

Here is first an example of  Theorem \ref{thmbande}.
{\ex \label{ex1} For $a=8$, $b=5$, we have $ab=40$. We then get:
\psfrag{x}{$\alpha$}
\psfrag{y}{$\beta$}
\begin{center}
\includegraphics[width=6cm]{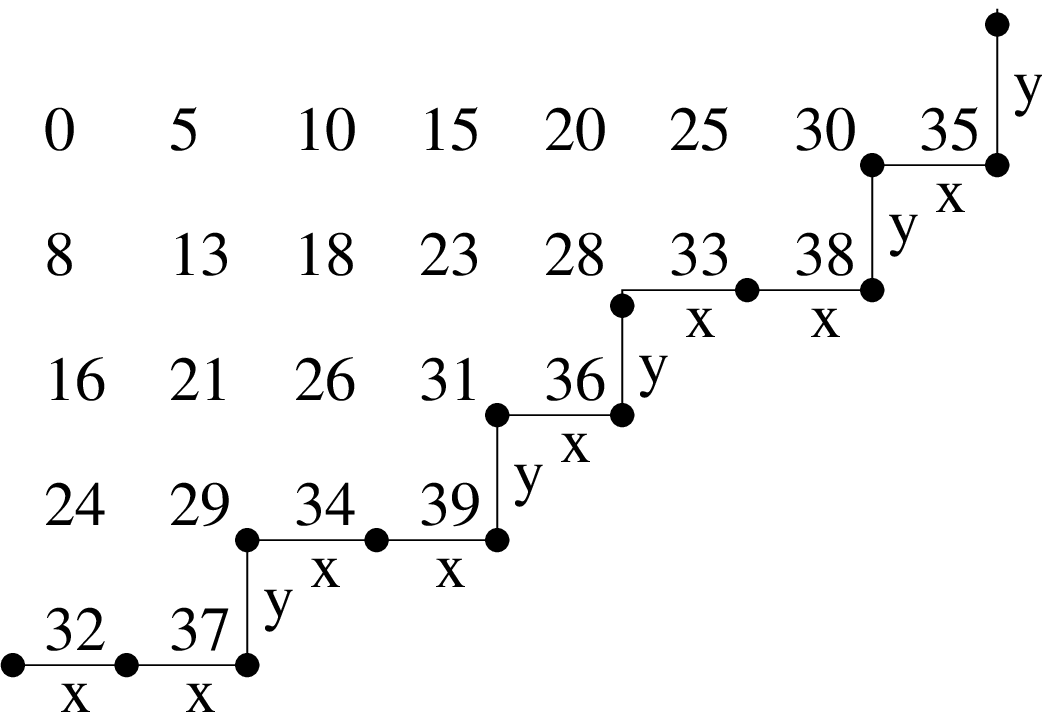}
\end{center}
Associating the letter $\alpha$ to a move to the right and the letter $\beta$ to a move to the top, and if we start at the lower leftmost corner, the lower boundary is coded by the word $\alpha \alpha \beta \alpha \alpha \beta \alpha \beta \alpha \alpha \beta \alpha \beta$: it is the Christoffel word with $8$ occurrences  of $\alpha$ and $5$ occurrences of $\beta$.
}

\Proof (of Theorem \ref{thmbande}) Let us consider the Cayley graph of the Christoffel word with $a$ occurrences of the letter $\alpha$ and $b$ occurrences of the letter $\beta$, with $\alpha < \beta$. We get the Cayley graph linearly represented by 
$$0 \rightarrow b \rightarrow 2b \,  \mod \,  (a+b) \rightarrow \ldots \rightarrow ib \, \mod \,  (a+b) \rightarrow \ldots \rightarrow (a+b-1)b \, \mod \,  (a+b) \rightarrow 0$$
In this Cayley graph, if there exists $k \in \N$ such that   
$$ib < k(a+b) \leq  (i+1)b,$$
 then 
\begin{equation} \label{eqma} (i+1)b \, \mod \,  (a+b) = (ib\, \mod \,  (a+b) ) - a.\end{equation}
Otherwise, we have 
\begin{equation}\label{eqpb} (i+1)b \, \mod \,  (a+b) = (ib \, \mod \,  (a+b))+b.\end{equation}
Let us consider the preceding Cayley graph to which we add the value $ab-a-b$. Since the values in the initial Cayley graph were lower or equal to $a+b$, the values in the new Cayley graph are now lower or equal to $a+b +ab-a-b=ab$. This corresponds exactly to take the lower and rightmost path such that the value of the coordinate $(x,-y)$ is lower or equal to $ab$. Indeed, we do $+b$ (see Equation (\ref{eqpb}): right move) if we exceed the value $ab$, otherwise we do $-a$ (see Equation (\ref{eqma}): up move). 
\QED

In the preceding example, the Cayley graph is
$$ 0 \rightarrow 5 \rightarrow 10 \rightarrow 2 \rightarrow 7 \rightarrow 12 \rightarrow 4 \rightarrow 9 \rightarrow 1 \rightarrow 6 \rightarrow 11 \rightarrow 3 \rightarrow 8 \rightarrow 0$$
The new Cayley graph obtained by adding $ab-a-b=27$ is 
$$ 27 \rightarrow 32 \rightarrow 37 \rightarrow 29 \rightarrow 34 \rightarrow 39 \rightarrow 31 \rightarrow 36 \rightarrow 28 \rightarrow 33 \rightarrow 38 \rightarrow 30 \rightarrow 35 \rightarrow 27$$
and corresponds to the boundary described in Example \ref{ex1}.

Here is a new proof of Corollary \ref{corMonnaie} that uses the result of Theorem \ref{thmbande}.

\Proof (of Corollary \ref{corMonnaie}) Excluding the integers that are on the boundary in Theorem \ref{thmbande} and using the Cayley graph seen previously, we obtain that there are exactly $xa+yb$ integers that are lower than $(a-1)(b-1)$. The total number of elements in the rectangle is $ab$ and since we have to remove the boundary which contains $a+b-1$ elements, and divide by  2, we obtain: $ \displaystyle \frac{ab-(a+b-1)}{2}=\frac{(a-1)(b-1)}{2}$.
\QED

\section{Concluding remarks}

In this paper, we have {expressed} in term of words combinatorics, a necessary and sufficient condition for the superimposition of two Christoffel words, by translating the results of \cite{rm1985,rjs2004} in terms of Christoffel words. For two superimposable Christoffel words,  we did more than in  \cite{rm1985,rjs2004} by giving a possible shift that always allows the perfect superimposition of two superimposable Christoffel words and  the number of possible shifts.  Those results are interesting since they give new properties of the well-known Christoffel words. Finally, in order to prove the Fraenkel conjecture, it would be interesting to generalize this result to the superimposition of more than two Christoffel words.

\bibliographystyle{alpha} 
{\footnotesize
\bibliography{soumission}
}

\end{document}